\documentclass[12pt]{article}
\usepackage{t1enc}
\usepackage[latin1]{inputenc}
\usepackage[english]{babel}
\usepackage{latexsym}
\newcommand{\bp}{{\bf P}}

\newcommand{\bs}{{\bf S}}

\newcommand{\bo}{{\bf 0}}

\newtheorem{lemma}{Lemma}[section]

\def\kx{{\bf x}}

\def\bo{{\bf 0}}

\def\bs{{\bf S}}

\def\zd{{\cal Z}_d}
\def\z2{{\cal Z}_2}

\def\begg{\begin{equation}}
\def\endd{\end{equation}}
\def\bege{\begin{eqnarray}}
\def\ende{\end{eqnarray}}

\def\ga{{\gamma}}

\begin{document}

\centerline{\Large\bf Heavy
points of a d-dimensional} \bigskip \centerline{\Large\bf simple random
walk}

\bigskip \bigskip \bigskip \bigskip \bigskip

\renewcommand{\thefootnote}{1} \noindent
\textbf{Endre Cs\'{a}ki}\footnote{Corresponding author. Research supported 
by the Hungarian National Foundation for Scientif\/ic Research, Grant No. 
T 037886 and T 043037.}\newline
Alfr\'ed R\'enyi Institute of Mathematics, Hungarian Academy of Sciences,
Budapest, P.O.B. 127, H-1364, Hungary. E-mail address: csaki@renyi.hu

\bigskip

\renewcommand{\thefootnote}{2} \noindent \textbf{Ant\'{o}nia
F\"{o}ldes}\footnote{Research supported by a PSC CUNY Grant,
No. 66494-0035.}\newline
Department of Mathematics, College of Staten Island, CUNY,
2800 Victory Blvd., Staten Island, New York
10314, U.S.A. E-mail address: foldes@mail.csi.cuny.edu

\bigskip

\renewcommand{\thefootnote}{3}
\noindent \textbf{P\'al R\'ev\'esz}\footnote{Research supported 
by the Hungarian National Foundation for Scientif\/ic Research, Grant No. 
T 037886 and T 043037.} 
\newline
Institut f\"ur Statistik und Wahrscheinlichkeitstheorie, Technische
Universit\"at Wien, Wiedner Hauptstrasse 8-10/107 A-1040 Vienna, Austria.
E-mail address: reveszp@renyi.hu

\bigskip \bigskip \bigskip

\noindent \textit{Abstract}: For a simple symmetric random walk in
dimension $d\geq 3$, a uniform strong law of large numbers is proved for
the number of sites with given local time up to time $n$.

\bigskip


\bigskip

 \vspace{.1cm}

\vfill
\newpage \renewcommand{\thesection}{\arabic{section}.}

\section{Introduction and main results}

\renewcommand{\thesection}{\arabic{section}} \setcounter{equation}{0}
\setcounter{theorem}{0} \setcounter{lemma}{0}

\noindent Consider a simple symmetric random walk
$\{\mathbf{S}_n\}_{n=1}^{\infty} $ starting at the origin $\mathbf{0}$
on the $d$-dimensional integer lattice $\mathcal{Z}_d$, i.e.
$\mathbf{S}_0=\mathbf{0}$, $\mathbf{S}_n=\sum_{k=1}^{n} \mathbf{X}_k$,
$n=1,2,\dots$, where $\mathbf{X}_k,\, k=1,2,\dots$ are i.i.d. random
variables with distribution
\[
\mathbf{P} (\mathbf{X}_1=\mathbf{e}_i)=\mathbf{P} (\mathbf{X}_1=
-\mathbf{e}_i)=\frac{1}{2d},\qquad i=1,2,...,d
\]
and $\{\mathbf{e}_1,\mathbf{e}_2,...\mathbf{e}_d\}$ is a system of
orthogonal unit vectors in $\mathcal{Z}_d.$  Define the local time of the
walk by
\begin{equation}
\xi(\mathbf{x},n):=\#\{k: \,\,0< k \leq n,\,\,\, \mathbf{S} _k=\mathbf{x}
\},\quad n=1,2,\ldots,
\label{loc1}
\end{equation}
where $\mathbf{x}$ is any lattice point of $\mathcal{Z}_d.$ The maximal
local time of the walk is defined as
\begin{equation}
\xi(n):=\max_{\mathbf{x} \in \mathcal{Z}_d}\xi(\mathbf{x},n).  \label{loc2}
\end{equation}

Define also
\begin{equation}
\eta(n):=\max_{0\leq k\leq n}\xi(\mathbf{S}_k,\infty).  \label{loc3}
\end{equation}

Denote by $\gamma(n)=\gamma(n;d)$ the probability that in the first $n-1$
steps the $d$-dimensional path does not return to the origin. Then
\begin{equation}
1=\gamma(1)\geq \gamma(2)\geq ...\geq \gamma(n)\geq...>0.
\label{gain}
\end{equation}
It was proved in \cite{DE50} that

\noindent
\textbf{Theorem A} (Dvoretzky and Erd\H os \cite{DE50}) \textit{For $d\geq 3
$
\begin{equation}
\lim_{n\to\infty}\gamma(n) =\gamma=\gamma(\infty;d)>0,  \label{gam1}
\end{equation}
and
\begin{equation}
{\gamma}<\gamma(n)<{\gamma}+O(n^{1-d/2}),  \label{gam2}
\end{equation}
or equivalently
\begin{equation}
\mathbf{P}(\xi(\mathbf{0},n)=0,\, \xi(\mathbf{0},\infty)>0)=O\left(
n^{1-d/2}\right)
\end{equation}
as} $n\to\infty$.

So $\gamma$ is the probability that the $d$-dimensional simple symmetric
random walk never returns to its starting point.

Let $\xi(\mathbf{x},\infty)$ be the total local time at $\mathbf{x}$ of the
infinite path in $\mathcal{Z}_d$. Then (see Erd\H os and Taylor
\cite{ET60}) $\,\xi(\mathbf{0},\infty)$ has geometric distribution:

\begin{equation}
\mathbf{P}(\xi(\mathbf{0},\infty)=k)={\gamma}(1-{\gamma})^k,\qquad
k=0,1,2,...  \label{geo1}
\end{equation}
Erd\H{o}s and Taylor \cite{ET60} proved the following strong law for the
maximal local time:

\bigskip \noindent \textbf{Theorem B} (Erd\H os and Taylor \cite{ET60})
\textit{For $d\ge 3$
\begin{equation}
\lim_{n\to\infty}\frac{\xi(n)}{\log n}=\lambda \hspace{1cm} \mathrm{a.s.},
\label{la}
\end{equation}
where}
\begin{equation}
\lambda=\lambda_d=-\frac{1}{\log(1-{\gamma})}.
\end{equation}

\bigskip Following the proof of Erd\H{o}s and Taylor, without any new
idea,
one can prove that
\begin{equation}
\lim_{n\to\infty}\frac{\eta(n)}{\log n}=\lambda \hspace{1cm} \mathrm{a.s.}
\label{lb}
\end{equation}

We can present a stronger lower estimate of $\xi(n)$.

\bigskip\noindent \textbf{Theorem C} (R\'ev\'esz \cite{R04}) \textit{Let $
d\geq 4$ and
\begin{equation}
\psi(n)=\psi(n,B)=\lambda\log n-\lambda B\log\log n.
\label{psi}
\end{equation}
Then with probability 1 for any $\varepsilon>0$ there is a random 
variable $n_0$ such that
\[
\xi(n)\geq\psi(n,3+\varepsilon)
\]
if $n\geq n_0$.}

\bigskip Erd\H{o}s and Taylor \cite{ET60} also investigated the properties
of
\[
Q(k,n):=\#\{\mathbf{x}:\ \mathbf{x}\in\mathcal{Z}_d,\
\xi(\mathbf{x},n)=k\},
\]
i.e. the cardinality of the set of points visited exactly $k$ times in the
time interval $[1,n]$. They proved

\bigskip\noindent \textbf{Theorem D} (Erd\H os and Taylor \cite{ET60})
\textit{
For $d\geq 3$ and for any} $k=1,2,\ldots$
\begin{equation}
\lim_{n\rightarrow\infty}{\frac{Q(k,n)}{n}}=
\gamma^2(1-\gamma)^{k-1}
\hspace{1cm} \mathrm{a.s.}  \label{qkn}
\end{equation}

Let
\begin{eqnarray}
U(k,n)&:=&
\#\{j:\ 0< j\leq n,\ \xi(\mathbf{S}_j,\infty)=k,\ \mathbf{S}_j\neq
\mathbf{S}_\ell\ (\ell=1,2,\ldots,j-1)\}\nonumber \\
&=&\#\{\kx\in\zd:\, 0<\xi(\kx,n)\leq \xi(\kx,\infty)=k\}.
\end{eqnarray}

Repeating the proof of Theorem D one can get
\begin{equation}
\lim_{n\rightarrow\infty}{\frac{U(k,n)}{n}}
=\gamma^2(1-\gamma)^{k-1}\hspace{1cm} \mathrm{a.s.}  \label{ukn}
\end{equation}
for any $k=1,2,\ldots$.

Define furthermore

\begin{eqnarray}
&&R(k,n):=\sum_{j=k}^\infty Q(j,n),\\
&&V(k,n):=\sum_{j=k}^\infty U(j,n).
\end{eqnarray}

It follows that for fixed $k\geq 1$
\begin{eqnarray}
&&\lim_{n\rightarrow\infty}{\frac{R(k,n)}{n}}
=\gamma(1-\gamma)^{k-1}\hspace{1cm} \mathrm{a.s.}  \label{rkn}\\
&&\lim_{n\rightarrow\infty}{\frac{V(k,n)}{n}}
=\gamma(1-\gamma)^{k-1}\hspace{1cm} \mathrm{a.s.}  \label{vkn}
\end{eqnarray}

The properties of these quantities were further investigated (for fixed
$k$) by Pitt \cite{Pitt} who proved (\ref{qkn}), (\ref{ukn}) and
(\ref{rkn}), (\ref{vkn}) for general random walk and by
Hamana \cite{Ham1}, \cite{Ham2} who proved central limit theorems (in
general case for $d\geq 3$).

In this paper we study the question whether $k$ can be replaced by a
sequence $t(n)=t_n\nearrow\infty$ of positive integers in (\ref{qkn}),
(\ref{ukn}), (\ref{rkn}) and (\ref{vkn}).

\bigskip\noindent
\textbf{Theorem}\textit{\,\, Let $d\geq 3$, and define
\begin{eqnarray}
\mu=\mu(t)&:=&\gamma(1-\gamma)^{t-1},\label{mut}\\
t_n&:=&[\psi(n,B)],\quad B>2,\label{tn}
\end{eqnarray}
where $\psi(n,B)$ is defined by {\rm (\ref{psi})}. Then we have}
\begin{eqnarray}
&&\lim_{n\rightarrow\infty}\sup_{t\leq
t_n}\left|{\frac{U(t,n)}{n\gamma\mu(t)}}-
1\right|=0\hspace{1cm} \mathrm{a.s.}  \label{utn}\\
&&\lim_{n\rightarrow\infty}\sup_{t\leq
t_n}\left|{\frac{Q(t,n)}{n\gamma\mu(t)}}-
1\right|=0\hspace{1cm} \mathrm{a.s.}  \label{qtn}\\
&&\lim_{n\rightarrow\infty}\sup_{t\leq t_n}\left|{\frac{V(t,n)}{n\mu(t)}}-
1\right|=0\hspace{1cm} \mathrm{a.s.}  \label{vtn}\\
&&\lim_{n\rightarrow\infty}\sup_{t\leq t_n}\left|{\frac{R(t,n)}{n\mu(t)}}-
1\right|=0\hspace{1cm} \mathrm{a.s.}  \label{rtn}
\end{eqnarray}
Here in $\sup_{t\leq t_n}$, $t$ runs through positive integers.

(\ref{rtn}) of Theorem  clearly implies (compare to Theorem C)

\bigskip\noindent
\textbf{Corollary}
\textit{\,\, Let $d\geq 3$. Then with probability 1 for any 
$\varepsilon>0$ there is a random variable $n_0$ such that
\[
\xi(n)\geq\lambda\log n-(2+\varepsilon)\log\log n
\]
if $n\geq n_0$.}

\bigskip
First we present some more notations. For $\mathbf{x}\in \mathcal{Z}_d$ let
$T_\kx$ be the first hitting time of $\mathbf{x}$, i.e.
$T_\kx=\min\{i\geq 1:\mathbf{S}_i=\mathbf{x}\}$ with the convention that
$T_\kx=\infty$ if there is no $i$ with $\mathbf{S}_i=\mathbf{x}$.
Let $T=T_\bo$. In general, for a subset $A$ of $\mathcal{Z}_d$,
let $T_A$ denote the first time the random
walk visits $A$, i.e. $T_A=\min\{i\geq 1:\, \mathbf{S}_i\in A\}
=\min_{\mathbf{x}\in A}T_{\mathbf{x}}$. Let
$\mathbf{P}_{\mathbf{x}}(\cdot)$ denote the probability of the
event in the bracket under the condition that the random walk starts from
$\mathbf{x}\in\mathcal{Z}_d$. We denote $\mathbf{P}(\cdot)=
\mathbf{P}_{\mathbf{0}}(\cdot)$.

Introduce further
\begin{eqnarray}
q_{\mathbf{x}}&:=&\mathbf{P}(T<T_{\mathbf{x}}),\label{aq} \\
s_{\mathbf{x}}&:=&\mathbf{P}(T_{\mathbf{x}}<T).\label{as}
\end{eqnarray}

In words, $q_{\mathbf{x}}$ is the probability that the random walk,
starting from $\mathbf{0}$, returns to $\mathbf{0}$, before reaching
$\mathbf{x}$ (including $T<T_{\mathbf{x}}=\infty$), and
$s_{\mathbf{x}}$ is the probability that the random walk, starting from
$\mathbf{0}$, hits $\mathbf{x}$, before returning to $\mathbf{0}$
(including $T_{\mathbf{x}}<T=\infty$).

\renewcommand{\thesection}{\arabic{section}.}

\section{Preliminary facts and results}

\renewcommand{\thesection}{\arabic{section}} \setcounter{equation}{0}
\setcounter{theorem}{0} \setcounter{lemma}{0}

First we present some lemmas needed to prove Theorem.

Introduce the following notations:

\newpage
\begin{eqnarray*}
X_i(t)&=&X_i= \\
&=&\left\{
\begin{array}{ll}
1 & \mathrm{if\ }\mathbf{S}_j\neq \mathbf{S}_i\ (j=1,2,\ldots,i-1),\
\xi(\mathbf{S}_i,\infty)\geq t, \\
0 & \mathrm{otherwise},
\end{array}
\right. \\
Y_i(t,n)&=&Y_i= \\
&=&\left\{
\begin{array}{ll}
1 & \mathrm{if\ } \mathbf{S}_j\neq \mathbf{S}_i\ (j=1,2,\ldots,i-1),\
\xi(\mathbf{S}_i,n)\geq t, \\
0 & \mathrm{otherwise},
\end{array}
\right. \\
\rho_i&=&\rho_i(t)=I\{X_i=1\}(\min\{j:\ \xi(\mathbf{S}_i,j)\geq t\}-i), \\
\mu_i&=&\mu_i(t)=\gamma(i)(1-\gamma)^{t-1}, \\
\end{eqnarray*}
$t=1,2,\ldots$, $i=1,2,\ldots$, where $I\{\cdot\}$ denotes the usual
indicator function.

Recall the definitions of $\gamma(i),$ $\gamma$ and $\mu=\mu(t)$ in (\ref
{gain}) (\ref{gam1}) and (\ref{mut}). Furthermore let
\begin{equation}
\sigma_n^2=\sigma_n^2(t):=\mathbf{E}\left(\sum_{i=1}^n X_i-n\mu\right)^2.
\end{equation}
Clearly we have
\[
R(t,n)=\sum_{i=1}^n Y_i,
\]
\[
V(t,n)=\sum_{i=1}^n X_i.
\]

\begin{lemma}
{\rm(Dvoretzky and Erd\H{o}s \cite{DE50})}

\[
\mathbf{P}(\mathbf{S}_i\neq \mathbf{S}_j,\ j=1,2,\ldots,i-1)=\mathbf{P}
(\xi(\mathbf{0},i-1)=0)=\gamma(i).
\]
\end{lemma}

The following lemma is a trivial consequence of Theorem A.

\begin{lemma}
\[
\mathbf{P}(n<\rho_i(t)<\infty)\leq{\frac{O(1)t^{d/2}}{n^{d/2-1}}},
\]
\[
\mu\leq\mu_i\leq\left(1+{\frac{O(1)}{i^{d/2-1}}}\right)\mu,
\]
\[
\mathbf{E} X_i=\mu_i.
\]
\end{lemma}

The next lemma can be obtained by elementary calculations.

\begin{lemma}
\[
n\mu\leq\mathbf{E}\sum_{i=1}^n X_i=\sum_{i=1}^n\mu_i\leq n\mu+\mu a_nO(1),
\]
where
\[
a_n=\sum_{i=1}^n{\frac{1}{i^{d/2-1}}}=\left\{
\begin{array}{ll}
O(1) & \mathrm{if}\quad d>4, \\
O(1)\log n & \mathrm{if}\quad d=4, \\
O(1)n^{1/2} & \mathrm{if}\quad d=3.
\end{array}
\right.
\]
\end{lemma}

\begin{lemma} Let $n>3^3$. Then
\begin{equation}
\sigma_n^2\leq n\mu+\mu a_nO(1)-n^2\mu^2+2(I+II+III),
\label{sig}
\end{equation}
where
\begin{eqnarray*}
I&=&\sum_{1\leq i<j\leq n}\mathbf{P}(X_i=1,\ X_j=1,\ \rho_i\geq
n^\alpha),
\\
II&=&\sum_{1\leq i<j\leq \min(i+3n^\alpha,n)}\mathbf{P}(X_i=1,\ X_j=1,\
\rho_i<n^\alpha), \\
III&=&\sum_{1\leq i<i+3n^\alpha<j\leq n}\mathbf{P}(X_i=1,\ X_j=1,\
\rho_i<n^\alpha), \\
\alpha&=&2/d.
\end{eqnarray*}
\end{lemma}

\vspace{2ex}\noindent \textbf{Proof.} Clearly we have
\begin{eqnarray*}
\sigma_n^2&=&\mathbf{E}\left(\sum_{i=1}^n
X_i\right)^2+n^2\mu^2-2n\mu\mathbf{E}
\sum_{i=1}^n X_i= \\
&=&\mathbf{E}\sum_{i=1}^n X_i+2\sum_{1\leq i<j\leq n}\mathbf{E}
X_iX_j+n^2\mu^2-2n\mu\sum_{i=1}^n\mu_i\leq \\
&\leq &n\mu+\mu a_nO(1)+2\sum_{1\leq i<j\leq n}\mathbf{E} X_iX_j-n^2\mu^2.
\end{eqnarray*}
Further
\[
\sum_{1\leq i<j\leq n}\mathbf{E} X_iX_j=\sum_{1\leq i<j\leq n} \mathbf{P}
\{X_i=1,\ X_j=1\}=I+II+III.
\]

Hence Lemma 2.4 is proved.

\medskip
Now let $A^{(\mathbf{x})}$ denote the two-point set
$\{\mathbf{0},\mathbf{x}\}$ and let $\Xi(A^{(\mathbf{x})},\infty)=
\xi(\mathbf{0},\infty)+\xi(\mathbf{x},\infty)$ denote its total occupation
time.

\begin{lemma} For $\kx\in\zd,\, \kx\neq \bo$, define
$\gamma_{\kx}:=\mathbf{P}(T_{\kx}=\infty)$ and recall the definitions of
$q_{\kx}$ and $s_{\kx}$ in {\rm (\ref{aq})} and {\rm (\ref{as})}. Then
\begin{eqnarray}
{\gamma}_{\mathbf{e}_i}&=&{\gamma}_{-\mathbf{e}_i}={\gamma},\quad
i=1,2,\ldots,d,  \label{eg} \\
\gamma_{\mathbf{x}}&\geq&\gamma, \label{xe} \\
q_{\mathbf{x}}&=&1-\frac{{\gamma}}{1-(1-{\gamma}_\kx)^2},  \label{qu} \\
s_{\mathbf{x}}&=&(1-{\gamma}_{\mathbf{x}})(1-q_\kx),  \label{es} \\
q_{\mathbf{x}}+s_{\mathbf{x}}&=&1-\frac{\ga}{2-\ga_{\kx}},\label{esqu} \\
\mathbf{P}(\Xi(A^{(\mathbf{x})},\infty)=j)&=&
(1-q_\kx-s_\kx)(q_\kx+s_\kx)^j,\quad
j=0,1,\ldots.  \label{dgeo5}
\end{eqnarray}
\end{lemma}

\vspace{2ex}\noindent \textbf{Proof.}
We show (\ref{eg}) first. For symmetric reason,
${\gamma}_{\pm\mathbf{e}_i}={\gamma}_{\pm\mathbf{e}_j}$, $i,j=1,\ldots,d$.
Hence
\[
1-{\gamma}=\sum_{i=1}^d \mathbf{P}(\mathbf{S}_1=\mathbf{e}_i)
(1-{\gamma}_{\mathbf{e}_i})+ \sum_{i=1}^d \mathbf{P}(\mathbf{S}_1=
-\mathbf{e}_i)(1-{\gamma}_{-\mathbf{e}_i})
=2\sum_{i=1}^d\frac1{2d}(1-{\gamma}_{\mathbf{e}_1})
=1-{\gamma}_{\mathbf{e}_1},
\]
proving (\ref{eg}).

To show (\ref{xe}), observe that starting from the origin, before hitting
$\mathbf{x}$ with $\Vert\mathbf{x}\Vert>1$, the random walk should hit
first the sphere $S(\mathbf{x},1):=\{\mathbf{y}:\,
\Vert\mathbf{y}-\mathbf{x}\Vert=1\}$. Hence
\begin{equation}
1-{\gamma}_{\mathbf{x}}=\mathbf{P}(T_{S(\mathbf{x},1)}<\infty)
(1-{\gamma})\leq 1-{\gamma}.  \label{fontos}
\end{equation}

Now let $Z(A)$ denote the number of visits in the set $A$ up to the first
return to zero, i.e.
\begin{equation}
Z(A)=\sum_{n=1}^T I\{\mathbf{S}_n\in A\}.
\end{equation}

Observe that

\begin{eqnarray}
\mathbf{P}(Z(A^{(\mathbf{x})})=j+1,T<\infty)= \left\{
\begin{array}{ll}
q_{\mathbf{x}} & \quad\mathrm{if \, \,}j=0 , \\
s_{\mathbf{x}}^2q_{\mathbf{x}}^{j-1} & \quad\mathrm{if \, \,} j=1,2,...
\label{gg2}
\end{array}
\right.
\end{eqnarray}
Summing up in (\ref{gg2}) we get
\begin{equation}
\sum_{j=0}^\infty\mathbf{P}(Z(A^{(\mathbf{x})})=j+1,T<\infty)=
q_{\mathbf{x}}+\frac{s_{\mathbf{x}}^2}
{1-q_{\mathbf{x}}}=\mathbf{P}(T<\infty)=1-\gamma.
\label{gg3}
\end{equation}
On the other hand, one can easily see that
\begin{eqnarray*}
1-\gamma&=&\mathbf{P}(T<\infty)=\mathbf{P}(T<T_{\mathbf{x}})+
\mathbf{P}(T>T_{\mathbf{x}},T< \infty) \\
&=&\mathbf{P}(T<T_{\mathbf{x}})+\mathbf{P}(T>T_{\mathbf{x}})
\mathbf{P}_{\mathbf{x}}(T<\infty) \\
&=&\mathbf{P}(T<T_{\mathbf{x}})+\mathbf{P}(T>T_{\mathbf{x}})
\mathbf{P}(T_{\mathbf{x}}
< \infty)=q_{\mathbf{x}}+ s_{\mathbf{x}}(1-\gamma_\kx),
\end{eqnarray*}
i.e.
\begin{equation}
1-\gamma=q_{\mathbf{x}}+s_{\mathbf{x}}(1-\gamma_\kx)  \label{gg4}
\end{equation}
Now (\ref{gg3}) and (\ref{gg4}) easily imply (\ref{qu}) and (\ref{es}),
hence also (\ref{esqu}).

Equation (\ref{dgeo5}) was proved in \cite{CFRRS} for general random walk.
For completeness a short proof is presented here. The probability that the
random walk, starting from $\mathbf{0}$, returns to $\mathbf{0}$ without
hitting $\mathbf{x}$, is $q_{\mathbf{x}}$, while $s_{\mathbf{x}}$ is the
probability that the random walk starting from $\mathbf{0}$ hits
$\mathbf{x}$ without returning to $\mathbf{0}$. Similarly, for
symmetric reason, $q_{\mathbf{x}}$ is also the probability of the
random walk starting from $\mathbf{x}$ returns to $\mathbf{x}$ without
hitting $\mathbf{0}$, and $s_{\mathbf{x}}$ is also the probability of the
random walk starting from $\mathbf{x}$ hits $\mathbf{0}$ in finite time,
without returning to $\mathbf{x}$. Hence, the probability that the random
walk starting from any point of $A^{(\mathbf{x})}$, returns to
$A^{(\mathbf{x})}$ in finite time, is
$q_{\mathbf{x}}+s_{\mathbf{x}}$. This gives (\ref{dgeo5}).

\medskip
Similarly to Theorem A, we prove

\begin{lemma}
\begin{eqnarray}
1-\gamma_\kx(n):=\mathbf{P}(T_\kx<n)&=&1-{\gamma}_\kx+\frac{O(1)}{n^{d/2-1}},
\label{gn} \\
q_\kx(n):=\mathbf{P}(T<\min(n,T_\kx))&=&q_\kx+\frac{O(1)}{n^{d/2-1}},
\label{qn} \\
s_\kx(n):=\mathbf{P}(T_\kx<\min(n,T))&=&s_\kx+\frac{O(1)}{n^{d/2-1}},
\label{sn}
\end{eqnarray}
and $O(1)$ is uniform in $\kx$.
\end{lemma}

\vspace{2ex}\noindent \textbf{Proof.} For the proof of (\ref{gn})
see Jain and Pruitt \cite{JP}.

To prove (\ref{qn}) and (\ref{sn}), observe that
\begin{eqnarray*}
q_\kx-q_\kx(n)&=&\mathbf{P}(T<T_\kx,\,n\leq T<\infty) \leq
\mathbf{P}(n\leq T<\infty)=\gamma(n)-\gamma, \\
s_\kx-s_\kx(n)&=&\mathbf{P}(T_\kx< T,\,n\leq T_\kx <\infty) \leq
\mathbf{P}(n\leq T_\kx<\infty)=\gamma_\kx(n)-\gamma_\kx.
\end{eqnarray*}

\begin{lemma}
Let $i<j$. Then for $t\geq 1$ integer we have
\begin{equation}
\mathbf{P}(X_i=1,X_j=1)\leq C\mu^2\left(1+\frac{t^{d/(d-2)}}{(j-i)^{d/2}}
\left(\frac2{2-\gamma}\right)^{2t}\right),
\end{equation}
where $C$ is a constant, independent of $i,j,t$ and $\mu=\mu(t)
=\gamma (1-\gamma)^{t-1}$.
\end{lemma}

\vspace{2ex}\noindent \textbf{Proof.}
Using (\ref{dgeo5}) of Lemma 2.5, we get
$$
\mathbf{P}(X_i=1,X_j=1)
$$
$$
\leq\sum_{\kx\in\zd}
\mathbf{P}(\bs_j-\bs_i=\kx,\xi(\bs_i,\infty)
-\xi(\bs_i,i)+\xi(\bs_j,\infty)-\xi(\bs_j,i)
\geq 2t-1)
$$
$$
=\sum_{\kx\in\zd}\mathbf{P}(\mathbf{S}_{j-i}=\kx)
\mathbf{P}(\Xi(A^{(\kx)},\infty)\geq 2t-1)
$$
$$
=\sum_{\kx\in\zd}\mathbf{P}(\bs_{j-i}=\kx)
(q_{\kx}+s_{\kx})^{2t-1}=\sum_{\kx\in\zd,\Vert \kx\Vert\leq R}+
\sum_{\kx\in\zd, \Vert \kx\Vert>R},
$$
where $R$ will be chosen later. For estimating the first sum, we use
 $\gamma_{\kx}\geq\gamma$ (cf. (\ref{xe}) of Lemma 2.5), hence by
(\ref{esqu})
$$
q_{\kx}+s_{\kx}=1-\frac{\gamma}{2-\gamma_{\kx}}\leq \frac{2(1-\gamma)}
{2-\gamma}.
$$
On the other hand
$$
\mathbf{P}(\bs_{j-i}=\kx)\leq \frac{C_1}{(j-i)^{d/2}},\qquad \kx\in\zd
$$
with some constant $C_1$, not depending on $\kx$ (cf. Spitzer \cite{S76},
page 72).

Since the cardinality of the set $\{\Vert \kx\Vert\leq R\}$ is a constant
multiple of $R^d$, we have
\begin{equation}
\sum_{\kx\in\zd, \Vert \kx\Vert\leq R}
\leq \frac{C_2R^d}{(j-i)^{d/2}}\left(\frac{2(1-\gamma)}{2-\gamma}
\right)^{2t}
\label{smallx}
\end{equation}
with some constant $C_2$.

For estimating the second sum, we use $1-\gamma_{\kx}\leq C_3R^{-d+2}$
for $\Vert \kx\Vert >R$ (cf. R\'ev\'esz \cite{R90}, page 241), hence
$$
q_{\kx}+s_{\kx}\leq 1-\gamma+C_4R^{-d+2}=(1-\gamma)
\left(1+\frac{C_4}{(1-\gamma)R^{d-2}}\right).
$$
Now choose $R=t^{1/(d-2)}$. Then
$$
(q_{\kx}+s_{\kx})^{2t-1}\leq C_5 (1-\gamma)^{2t}.
$$
Here the constant $C_5$ is independent of both $\kx$ and $t$. Since
$$
\sum_{\kx\in\zd}\mathbf{P}(\bs_j-\bs_i=\kx)=1,
$$
we have
$$
\sum_{\kx\in\zd,\ \Vert \kx\Vert >R}\leq C_5(1-\gamma)^{2t}=C_6\mu^2.
$$
this together with (\ref{smallx}) (putting $R=t^{1/(d-2)}$ there) proves
Lemma 2.7.

In the subsequent lemmas $t_n$ is defined by (\ref{tn}).

\begin{lemma} For $t\leq t_n$, any $\varepsilon>0$ and large enough $n$ we
have
\begin{equation}
I\leq
O(1)n^{2/d+\varepsilon}\left(n+\left(\frac{2}{2-\gamma}\right)^{2t_n}\right)
\mu^2(t).
\label{one}
\end{equation}
\end{lemma}

\vspace{2ex}\noindent \textbf{Proof.}
Now we need to estimate the probability
\[
\mathbf{P}(X_i=1,X_j=1, \rho_i\geq n^\alpha).
\]

Define the events $B_k$ by
\[
B_k=\{\xi(\bs_i,\infty)-\xi(\bs_i,i)+\xi(\bs_j,\infty)-\xi(\bs_j,i)=k\}
\]
and consider the $k$ time intervals between the consecutive visits of
$\{\bs_i,\bs_j\}$. Then at least one of these intervals is larger than
\begin{equation}
{\frac{\rho_i(t)}{k}}\geq{\frac{n^\alpha}{k}}  \label{nex}
\end{equation}
(provided that $\{X_i=1,\ X_j=1,\ \rho_i\geq n^\alpha\}$). Denote this
event by $D_k$. Similarly to the proof of Lemma 2.7 we have
\[
\mathbf{P}(X_i=1,X_j=1, \rho_i\geq n^\alpha)
\leq \sum_{\kx\in\zd}\mathbf{P}(\bs_j-\bs_i=\kx,\, \cup_{k\geq
2t-1}B_kD_k)
\]
\[
\leq \sum_{\kx\in\zd}\mathbf{P}(\bs_{j-i}=\kx)
\sum_{k\geq 2t-1}\mathbf{P}(B_kD_k\, |\, \bs_j-\bs_i=\kx).
\]

The event $B_kD_k$, under the condition $\bs_j-\bs_i=\kx$,
means that placing a new origin at the point
$\bs_i$, and starting the time at $i$, there are exactly $k$ visits
in the set $A^{(\kx)}$, and at least one time interval between
consecutive visits is larger than $n^\alpha/k$. Hence applying
(\ref{dgeo5}) of Lemma 2.5 and (\ref{qn}), (\ref{sn}) of Lemma 2.6,
we get
\[
\mathbf{P}(B_kD_k\, |\, \bs_j-\bs_i=\kx)
\leq k(1-q_\kx-s_\kx)(q_\kx+s_\kx)^{k-1}
\left(q_\kx+s_\kx-q_\kx\left(\frac{n^\alpha}{k}\right)
-s_\kx\left(\frac{n^\alpha}{k}\right)\right)
\]
\[
\leq O(1)k\left(\frac{k}{n^\alpha}\right)^{d/2-1}
(1-q_\kx-s_\kx)(q_\kx+s_\kx)^{k-1}\leq
O(1)k^{d/2}n^{2/d-1}(q_\kx+s_\kx)^{k-1},
\]
where $O(1)$ is uniform in $k$ and $\kx$, hence
\begin{eqnarray*}
\sum_{k\geq 2t-1}\mathbf{P}(B_kD_k\, |\, \bs_j-\bs_i=\kx)
&\leq& O(1)n^{2/d-1}
\sum_{k\geq 2t-1}k^{d/2}(q_\kx+s_\kx)^{k-1}\\
&\leq& O(1)n^{2/d-1}t^{d/2}(q_\kx+s_\kx)^{2t-2}.
\end{eqnarray*}

Proceeding now as in the proof of Lemma 2.7, we can estimate
\[
\mathbf{P}(X_i=1,X_j=1, \rho_i\geq n^\alpha)\leq
O(1)t^{d/2}n^{2/d-1}\mu^2(t)\left(1+\frac{t^{d/(d-2)}}{(j-i)^{d/2}}
\left(\frac2{2-\gamma}\right)^{2t}\right)
\]
and summing up for $1\leq i<j\leq n$, we get
\[
I\leq
O(1)n^{2/d}t_n^{d/2}
\left(n+t_n^{d/(d-2)}\left(\frac{2}{2-\gamma}\right)^{2t_n}\right)
\mu^2(t),
\]
since $t\leq t_n$. But $t_n<\lambda\log n$, therefore any power of $t_n$
can be estimated by $n^\varepsilon$, hence (\ref{one}) follows.

\begin{lemma} For $t\leq t_n$, any $\varepsilon>0$ and large enough $n$ we
have
\begin{equation}
II\leq
O(1)n^{2/d+\varepsilon}\left(n+
n^{1-2/d}\left(\frac{2}{2-\gamma}\right)^{2t_n}\right)
\mu^2(t).
\label{two}
\end{equation}
\end{lemma}

\vspace{2ex}\noindent \textbf{Proof}.
Using the estimate in Lemma 2.7 and
summing up for $i,j$ with $1\leq i<j\leq \min(i+3n^\alpha,n)$,
using again that $t_n<\lambda\log n$, a simple calculation shows
(\ref{two}).

\begin{lemma} For $t\leq t_n$, any $\varepsilon>0$ and large enough $n$ we
have
\begin{equation}
III\leq\frac{\mu^2(t) n^2}{2}+O(1)n^{3/2}\mu^2(t).
\label{three}
\end{equation}
\end{lemma}

\vspace{2ex}\noindent \textbf{Proof.} Let
\begin{eqnarray*}
A&=&\{\mathbf{S}_i\ \mathrm{is\ a\ new\ point\ i.e.\ } \mathbf{S}_i\neq
\mathbf{S}_j\ j=1,2,\ldots,i-1\}, \\
B&=&\{\xi(\mathbf{S}_i,i+n^\alpha)-\xi(\mathbf{S}_i,i)\geq t-1\}, \\
D&=&\{\mathbf{S}_j\ \mathrm{is\ a\ new\ point}\}, \\
E&=&\{\xi(\mathbf{S}_j,\infty)-\xi(\mathbf{S}_j,j)\geq t-1\}, \\
D\subset G&=&\left\{\xi(\mathbf{S}_j,j)
-\xi\left(\mathbf{S}_j,i+{\frac{2(j-i)}{3}}\right)= 0\right\}, \\
B\subset H&=&\{\xi(\mathbf{S}_i,\infty)-\xi(\mathbf{S}_i,i)\geq t-1\}.
\end{eqnarray*}

Recall the definition of $\gamma(n)$ in Section 1 and
let $j>i+3n^\alpha$. Then
\begin{eqnarray*}
\lefteqn{ \bp\{X_i=1,\ X_j=1,\ \rho_i<n^\alpha\}\leq\bp\{ABDE\}\leq} \\
&\leq&\mathbf{P}(ABGE)=\mathbf{P}(A)\mathbf{P}(B)
\mathbf{P}(G)\mathbf{P}(E)\leq \\
&\leq&\mathbf{P}(A)\mathbf{P}(H)\mathbf{P}(G)\mathbf{P}(E)= \\
&=&\gamma(i+1)(1-\gamma)^{t-1}\gamma((j-i)/3)
(1-\gamma)^{t-1}.
\end{eqnarray*}

Clearly we have
\begin{eqnarray*}
III&\leq&\sum\gamma(i+1)(1-\gamma)^{t-1}
\gamma((j-i)/3)(1-\gamma)^{t-1}\leq \\
&\leq&\gamma^2(1-\gamma)^{2t-2} \sum\left(1+{\frac{O(1)}{(j-i)^{d/2-1}}}
\right) \left(1+{\frac{O(1)}{i^{d/2-1}}}\right)\leq \\
&\leq&\gamma^2(1-\gamma)^{2t-2}\left[{n \choose
2}+O(1)(K+L+M)\right]
\end{eqnarray*}
where the summations above and below go for $\{i,j:\, 1\leq 
i<i+3n^\alpha<j\leq n\}$ and
\begin{eqnarray*}
K&=&\sum{\frac{1}{i^{d/2-1}}}\leq na_n, \\
L&=&\sum{\frac{1}{(j-i)^{d/2-1}}}\leq na_n, \\
M&=&\sum{\frac{1}{i^{d/2-1}}}{\frac{1}{(j-i)^{d/2-1}}}\leq na_n.
\end{eqnarray*}
Using $a_n=O(1)n^{1/2}$ (see Lemma 2.3) we have (\ref{three}).

\begin{lemma} For $t\leq t_n$, any $\varepsilon>0$ and large enough $n$ we
have
\begin{equation}
\sigma_n^2= O(1)[n\mu(t)+\mu^2(t) n^{1.8}].
\label{sigma}
\end{equation}
\end{lemma}

\vspace{2ex}\noindent \textbf{Proof} is based on Lemmas 2.4,
2.8, 2.9 and 2.10. The numerical values of $\lambda$ can be obtained
by a result of Grif\/f\/in \cite{G90}:
\begin{eqnarray*}
1-\gamma_3&=&0.341, \\
1-\gamma_4&=&0.193, \\
1-\gamma_5&=&0.131, \\
1-\gamma_6&=&0.104.
\end{eqnarray*}
Consequently
\begin{eqnarray*}
\lambda_3&=&0.929, \\
\lambda_4&=&0.608, \\
\lambda_5&=&0.492, \\
\lambda_6&=&0.442.
\end{eqnarray*}

By using $t_n<\lambda \log n$, one can verify (numerically)
$$
\left(\frac2{2-\gamma}\right)^{2t_n}<
n^{2\lambda\log(2/(2-\gamma))}<n^{0.75}
$$
for $d=3$ and hence also for all $d\geq 3$. By choosing an appropriate
$\varepsilon$ and putting the estimations (\ref{one}), (\ref{two}),
(\ref{three}) into (\ref{sig}), we can see, that the term $n^2\mu^2$
cancels out and all the other terms are smaller than the right hand
side of (\ref{sigma}), proving Lemma 2.11.

Lemma 2.11 implies

\begin{lemma}
For any $0<C<B$, $t\leq t_n$ and large enough $n$ we have
\[
\sigma_n(\log n)^{C/2}\leq O(1)((n\mu(t))^{1/2}(\log n)^{C/2} +\mu(t)
n^{0.9}(\log n)^{C/2})=o(1)n\mu(t).
\]
\end{lemma}

 \renewcommand{\thesection}{\arabic{section}.}

\section{Proof of the Theorem }

\renewcommand{\thesection}{\arabic{section}} \setcounter{equation}{0}
\setcounter{theorem}{0} \setcounter{lemma}{0}

First we prove (\ref{vtn}).

By Markov's inequality for any $C>0$ we have
\[
\mathbf{P}(|V(t,n)-n\mu(t)|\geq\sigma_n(\log n)^{C/2})\leq(\log
n)^{-C}.
\]
By Lemma 2.12, if $C<B$,
\[
\mathbf{P}(|V(t,n)-n\mu(t)|\geq o(1)n\mu(t))\leq(\log n)^{-C}.
\]
Consequently, since $t_n<\lambda\log n$,
\begin{equation}
\mathbf{P}\left(\sup_{t\leq t_n+1}\frac{|V(t,n)-n\mu(t)|}{n\mu(t)}\geq
o(1)\right)\leq O(1)(\log n)^{-C+1}.
\label{supv}
\end{equation}

Choose $C>2,\ n(k)=\exp(k/\log k)$. (\ref{supv}) and Borel-Cantelli lemma
imply
\begin{equation}
\lim_{k\rightarrow\infty}\sup_{t\leq t(n(k))+1}
\left|{\frac{V(t,n(k))}{n(k)\mu(t)}}-1\right|=0\quad \mathrm{a.s.}
\label{supvtnk}
\end{equation}

Let $n(k)\leq n< n(k+1)$. Then for $t\leq t_n$ we have
\[
V(t,n(k))\leq V(t,n)\leq V(t,n(k+1))
\]
and
\[
\lim_{k\rightarrow\infty}{\frac{n(k+1)}{n(k)}}=1.
\]
Hence for any $\varepsilon>0$ and large enough $n$,
\[
\frac{V(t,n)}{n\mu(t)}\leq\frac{V(t,n(k+1))}{n(k+1)\mu(t)}
\frac{n(k+1)}n\leq (1+\varepsilon)\quad\mathrm{a.s.},
\]
since $t\leq t_n\leq t(n(k+1))$. Similarly,
\[
\frac{V(t,n)}{n\mu(t)}\geq\frac{V(t,n(k))}{n(k)\mu(t)} \frac{n(k)}n\geq
(1-\varepsilon)\quad\mathrm{a.s.}
\]
Hence we have (\ref{vtn}).

\medskip
Now we turn to the proof of (\ref{rtn}).

Let
$$
M(t,n)=V(t,n)-R(t,n)=\sum_{i=1}^n (X_i-Y_i).
$$
Observe that $X_i\geq Y_i$ and hence $M(t,n)$ is non-negative and
non-decreasing in $n$. Moreover, by Lemma 2.2
$$
\mathbf{E}(X_i-Y_i)=\mathbf{P}(X_i-Y_i=1)
\leq \mathbf{P}(X_i=1,n-i\leq\rho_i(t)<\infty)
\leq \frac{O(1)\mu(t) t^{d/2}}{(n-i)^{d/2-1}}.
$$
Consequently
$$
0\leq \frac{\mathbf{E}M(t,n)}{n\mu(t)}\leq\frac{O(1)(\log n)^{d/2}}
{n^{1/2}}.
$$

By Markov's inequality
$$
\mathbf{P}\left(\sup_{t\leq t_n}\frac{M(t,n)}{n\mu(t)}>\varepsilon\right)
\leq\frac{O(1)(\log n)^{d/2+1}}{n^{1/2}}.
$$
On choosing $n_k=k^{2+\delta}$, $\delta>0$, Borel-Cantelli lemma implies
$$
\lim_{k\to\infty}\sup_{t\leq t_{n_k}}\frac{M(t,n_k)}{n_k\mu(t)}=0
\qquad {\rm a.s.}
$$
Using the monotonicity of $M(t,n)$ in $n$, interpolating between
$n_k$ and $n_{k+1}$ we get
$$
\lim_{n\to\infty}\sup_{t\leq t_n}\frac{M(t,n)}{n\mu(t)}=0
\qquad {\rm a.s.}
$$
This combined with (\ref{vtn}) gives (\ref{rtn}).

(\ref{qtn}) and (\ref{utn}) are immediate from (\ref{rtn}) and
(\ref{vtn}), since $Q(t,n)=R(t,n)-R(t+1,n)$ and
$U(t,n)=V(t,n)-V(t+1,n)$.

This completes the proof of  the Theorem.


\begin{thebibliography}{99}
\bibitem{CFRRS} Cs\'aki, E., F\"oldes, A., R\'ev\'esz, P., Rosen, J.,
Shi, Z., 2005. Frequently visited sets for random walks. Stochastic 
Process. Appl., to appear.

\bibitem{DE50}  Dvoretzky, A., Erd\H os, P., 1951. Some problems on 
random walk in space. Proceedings of the Second Berkeley Symposium 
on Mathematical Statistics and Probability, University of California 
Press, Berkeley, pp. 353--367.

\bibitem{ET60}  Erd\H os, P., Taylor, S.J., 1960. Some problems 
concerning the structure of random walk paths. Acta Math. Acad. 
Sci. Hungar. 11, 137--162.

\bibitem{G90}  Griffin, P., 1990. Accelerating beyond the third 
dimension: returning to the origin in simple random walk. Math. 
Scientist 15, 24--35.

\bibitem{Ham1} Hamana, Y., 1992. On the central limit theorem for the 
multiple point range of random walk. J. Fac. Sci. Univ. Tokyo 39, 
339--363.

\bibitem{Ham2} Hamana, Y., 1995. On the multiple point range of three
dimensional random walk. Kobe J. 12, 95--122.

\bibitem{JP} Jain, N.C., Pruitt, W.E., 1971. The range of transient 
random walk. J. Analyse Math. 24, 369--393.

\bibitem{Pitt} Pitt, J.H., 1974. Multiple points of transient random 
walk. Proc. Amer. Math. Soc. 43, 195--199.

\bibitem{R90}  R\'ev\'esz, P., 1990. Random Walk in Random and Non-Random
Environments. World Scientific, Singapore.

\bibitem{R04}  R\'ev\'esz, P., 2004. The maximum of the local time of a 
transient random walk. Studia Sci. Math. Hungar. 41, 379--390.

\bibitem{S76}  Spitzer, F., 1976. Principles of Random Walk, 2nd. ed. Van
Nostrand, Princeton.
\end{thebibliography}
\end{document}